\documentclass[letterpaper, 10 pt, conference]{ieeeconf}
\IEEEoverridecommandlockouts
\usepackage{cite}
\usepackage{amsmath,amssymb,amsfonts}
\usepackage{algorithmic}
\usepackage{graphicx}
\usepackage{textcomp}
\usepackage{xcolor}

\usepackage{amsfonts, amsmath, amssymb, amsthm}
\def\BibTeX{{\rm B\kern-.05em{\sc i\kern-.025em b}\kern-.08em
    T\kern-.1667em\lower.7ex\hbox{E}\kern-.125emX}}

\usepackage{cleveref}

\usepackage{mathtools}
\newcommand{\overarc}[1]{\overset{\hbox{\large$\frown$}}{\hbox{ \vphantom{e}}}\mathllap{#1}}
\newcommand{\overarcSmall}[1]{\overset{\hbox{\scriptsize$\frown$}}{\hbox{ \vphantom{e}}}\mathllap{#1}}

\usepackage{alphabeta}
\newtheorem{theorem}{Theorem}
\newtheorem{lemma}{Lemma}
\newtheorem{proposition}{Proposition}
\theoremstyle{remark}
\newtheorem{remark}{Remark}
\usepackage[caption=false, font=footnotesize]{subfig}

\usepackage{float}

\begin{document}

\title{Pursuit-Evasion on a Sphere and When It Can Be Considered Flat}
\author{Dejan Milutinovi\'c, Alexander Von Moll, Satyanarayana G. Manyam, \\ David W. Casbeer, Isaac E. Weintraub, and Meir Pachter% <-this % stops a space
\thanks{
This paper is based on work performed at the Air Force Research Laboratory (AFRL) \textit{Control Science Center} and is supported by AFOSR LRIR 24RQCOR002.
        DISTRIBUTION STATEMENT A.
        Approved for public release.
        Distribution is unlimited.
        AFRL-2024-1467; 18 MAR 2024.}% <-this % stops a space
\thanks{D. Milutinovi\'c is with Electrical and Computer Engineering Dept., University of California, Santa Cruz, Santa Cruz, CA, USA {\tt\footnotesize dejan@ucsc.edu}}%
\thanks{A. Von Moll, D. Casbeer, and I. Weintraub are with the Control Science Center, Aerospace Systems Directorate, Air Force Research Laboratory, Dayton, OH 45435, USA {\tt\scriptsize alexander.von\_moll@us.af.mil, david.casbeer@us.af.mil, isaac.weintraub.1@us.af.mil}}%
\thanks{S. Manyam is with DCS Corporation, Dayton, OH, USA, {\tt\footnotesize smanyam@infoscitex.com}}%
\thanks{M. Pachter is with Dept. of Electrical Engineering, Air Force Institute of Technology, Wright-Patterson AFB, OH, USA, {\tt\footnotesize meir.pachter@afit.edu}}
}

\maketitle

\begin{abstract}
In classical works on a planar differential pursuit-evasion game with a faster pursuer, the intercept point resulting from the equilibrium strategies lies on the Apollonius circle. This property was exploited for the construction of the equilibrium strategies for two faster pursuers against one evader. Extensions for planar multiple-pursuer single-evader scenarios have been considered. We study a pursuit-evasion game on a sphere and the relation of the equilibrium intercept point to the Apollonius domain on the sphere. The domain is a generalization of the planar Apollonius circle set. We find a condition resulting in the intercept point belonging to the Apollonius domain, which is the characteristic of the planar game solution. Finally, we use this characteristic to discuss pursuit and evasion strategies in the context of two pursuers and a single slower evader on the sphere and illustrate it using numerical simulations. 
\end{abstract}

% \begin{IEEEkeywords}
% pursuit-evasion, sphere, Apollonius circle, non-Euclidean geometry,
% \end{IEEEkeywords}

\section{Introduction and Background}

\emph{The Apollonius circle} is a concept that has been used to 
address solutions of multiple non-cooperative pursuit-evasion (P-E) differential games under simple motion in the plane~\cite{Isaacs1965}. 
It is the set of points that can be reached by $E$ sooner than $P$, 
i.e., it is \emph{the dominance region of $E$}. In the same reference \cite{Isaacs1965}, the solution of the two pursuer and single slower evader differential game is inferred from the intersection of Apollonius circles. 
The problem was revisited and solved in \cite{Meir2019a} and the Apollonius circle construct was applied as a building block for the solution of a game with multiple pursuers \cite{Alex2020} in the plane. 
For P-E differential games on a sphere with a slower $E$, the dominance of $E$ is not a circle \cite{Kovshov2000}, therefore, we denote it as \emph{the Apollonius domain}.
\begin{figure}[htbp]
\begin{center}
\includegraphics[width=0.30\textwidth]{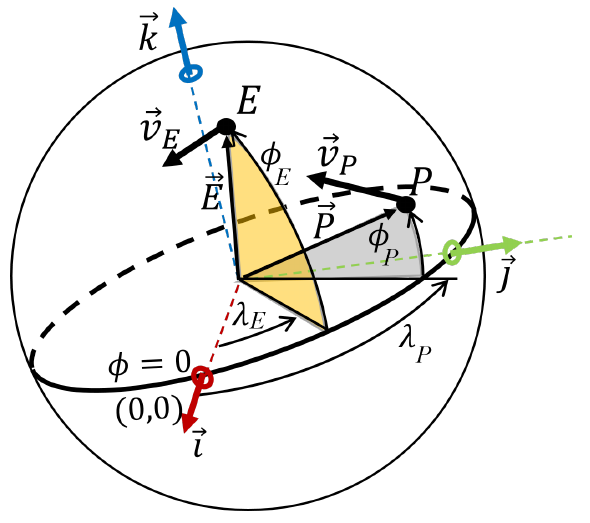}
\end{center}
\caption{Pursuit-evasion (P-E) of a pursuer (P) and a slower evader (E) on a sphere.}
\label{Fig1}
\end{figure}

This paper investigates the P-E game of $\min$-$\max$ capture time equilibrium strategies on a sphere for a slower $E$ and simple motion of both players, and examines their outcome against the Apollonius domain. 
In the planar case, the intercept point resulting from equilibrium strategies lies on the Apollonius circle, and here we investigate conditions for the P-E game on a sphere that result in the intercept point lying on the Apollonius domain. 
Under these conditions, we discuss P-E strategies for games on a sphere with multiple $P$s and a single $E$ in a similar fashion as in \cite{Alex2020}.

\emph{Related works}: If two players move in a plane and when both of them apply their equilibrium strategies, the corresponding equilibrium trajectories are straight lines. 
As pointed out in \cite{Melikyan2007} when a pursuer and an evader have simple motion over a 2D Riemannian manifolds, the same holds when the players are sufficiently close, except that instead of straight lines, the trajectories are along manifold shortest paths, i.e., geodesics. 
However, when the players are far from each other, there may be multiple geodesics of the same length, which requires an additional analysis. 
Using this perspective, the authors in \cite{Melikyan2007, Hovakimyan2000} analyzed multiple manifolds, but not a sphere. 
Note that when $P$ and $E$ are on the opposite sides of a sphere and on the same line with the center of the sphere, there are infinitely many geodesics between $P$ and $E$. 
The P-E game on a sphere was the topic of \cite{Kovshov2000}, but under the condition that $P$ always knows $E$'s velocity. Another related paper is \cite{Kuchkarov2009} which discusses sufficient conditions for a P-E solution. 
A P-E on a sphere also appears in  popular culture \cite{DBZ} in a graphic novel, in which, on a very small planet, the main character, Goku, must catch a pet monkey, Bubbles, to prove he is worthy to be trained by Master Kai. 

Since the P-E solution in this paper results in feedback strategies for $P$ and $E$, it is worth pointing to other works on feedback control on a sphere. 
These include the stabilization of collective motion \cite{Paley2009}, motion coordination in a spatio-temporal flow-field \cite{Paley2010}, the control of multi agents on a sphere to achieve different spherical patterns \cite{Spong2014} and consensus on a sphere \cite{Markdahl2018}, the stabilization of bilinear systems \cite{Saradagi2019} and collision-free motion coordination using barrier functions \cite{Egerstedt2020}. 

\emph{Paper contributions}: We revisit here the problem of P-E on a sphere, in which $P$ and a slower $E$ have simple motion. 
Our analysis relies on the Hamilton-Jacobi-Isaacs equation (HJI), which allows us to derive equilibrium strategies for both players, except for a special game configuration in which there are infinitely many geodesics connecting $P$ and $E$. 
The special configuration is the one in which $P$ and $E$ are on the opposite sides of the sphere and on the same line with the center of the sphere. 
Due to the multiplicity of solutions, the condition of the game is similar to the presence of the dilemma. 
Therefore, we use the value function and associated rate of loss \cite{TAC2023}, which allows us to obtain unique equilibrium strategies on the whole sphere surface. 
With this, we also obtain the corresponding intercept point and using results from \cite{Kovshov2000}, we are able to obtain the condition resulting in the intercept point in the Apollonius domain.

The paper is organized as follows. 
The problem formulation in Section \ref{sec:problem_description} is followed by the P-E game analysis and equilibrium strategies in Section \ref{sec:strategies}. 
Section  \ref{sec:Apollonius} is about a relation between the intercept point resulting from the equilibrium strategies and the Apollonius domain. 
Section \ref{sec:Results} shows examples illustrating paper results and Section \ref{sec:Conclusion} gives conclusions.

\section{Problem Description}
\label{sec:problem_description}
Consider two agents restricted to motion on a spherical surface, see Fig.~\ref{Fig1}: a pursuer ($P$) and a slower evader 
($E$). Both $P$ and $E$ move on the spherical surface and their positions are defined as the corresponding vectors
\begin{equation}
 \begin{aligned}
 \vec{P} &= R \cos \varphi_P \cos \theta_P \vec{i}+ R \cos \varphi_P \sin \theta_P \vec{j}+ R \sin \varphi_P\vec{k}, \\
 \vec{E} &= R \cos \varphi_E \cos \theta_E \vec{i}+ R \cos \varphi_E \sin \theta_E \vec{j}+ R \sin \varphi_E\vec{k}. 
 \end{aligned}
 \end{equation}

Due to their constrained motion on the sphere's surface, the vectors $\vec{E}$ and $\vec{P}$ defining the positions of $E$ and 
$P$ always have the same magnitudes $|\vec{E}|=|\vec{P}|=R$ (see Fig. \ref{Fig1}). The velocities $\vec{v}_E$ and $\vec{v}_P$ of $E$ and $P$, respectively, are 
\begin{align}
 \vec{v}_E &= v_E \vec{e},\ \  &0 \le v_E \le \mu v_P  \label{eq:veE}\\
 \vec{v}_P &= v_P \vec{p},\ \  &v_P=\mathrm{const},\ \ 0< \mu < 1,        \label{eq:vpP}
\end{align}
where $\vec{e}$ and $\vec{p}$ are unit vectors tangential to the sphere representing the agent's direction of travel, 
and $v_E$ and $v_P$ representing the corresponding speed; this implies $\vec{v}_E \cdot \vec{E}=\vec{v}_P \cdot \vec{P} = 0$.

The control variable of $P$ is only the velocity direction $\vec{p}$, while the control variables of $E$ are both the speed $v_E$ and the velocity direction $\vec{e}$. The speed ratio $\mu \in (0,1)$ is the ratio of $E$'s maximal speed and $P$'s speed $v_P$; therefore, $P$ is faster than $E$. 

\emph{P-E game}: The two agents are initially at different locations of the sphere and when \emph{$P$ is collocated with $E$}, we say that $P$ has \emph{captured} $E$. The capture time $\tau$ is the length of time that it takes $P$ and $E$ to move from their initial positions to when $P$ captures $E$, or formally
 \begin{equation}
   \tau = \inf \lbrace t \;| \; \vec{E}=\vec{P} \rbrace.
   \label{eq:tau}
 \end{equation}
The goal of $P$ is to select its sequence of actions to minimize the capture time, while the goal of $E$ is to select its sequence of actions that will maximize the capture time.

\emph{Problem}: We analyze here the pursuit evasion problem on the sphere and find 
the equilibrium (i.e., the saddle point or Nash-equilibrium) strategies 
$\{\vec{v}_P^*, \vec{v}_E^*\}$ for $P$ and $E$, respectively. They are 
optimal in a sense that if both $P$ and $E$ commit to them, then the 
capture time is $\tau^*(\vec{v}_P^*, \vec{v}_E^*)$. However, if $P$ ($E$) 
gives up its equilibrium strategy and applies any other action at any time point, 
it will extend (shorten) time to capture. This can be expressed as
\begin{equation}
\tau(\vec{v}_P^*, \vec{v}_E) \le \tau^*(\vec{v}_P^*, \vec{v}_E^*) \le 
\tau(\vec{v}_P, \vec{v}_E^*).
\end{equation}
The capture time $\tau^*$ is also called the Value of the game. 
Our goal is to find equilibrium strategies for $P$ and $E$, as well as the Value function. Then we will examine the relative position 
between the equilibrium intercept point (i.e., the intercept point resulting from both players applying their equilibrium strategies) and the Apollonius domain, which is defined as follows.

\emph{Definition 1}: \emph{Apollonius domain} $\mathcal{A}$. 
Let us assume that at the initial time $t=0$, the positions of $P$ 
and $E$ on a sphere are $P(0)$ and $E(0)$, $P(0)\ne E(0)$. 
The Apollonius domain, $\mathcal{A}$, is a closed set of points 
$I$ on the sphere defined as
\begin{equation}
 \mathcal{A} = \left\{ I \ | \ \tfrac{\overarcSmall{IE}(0)}{\mu v_P} \le \tfrac{\overarcSmall{IP}(0)}{v_P} \right\},
\end{equation}
where $\overarc{IE}(0)$, $\overarc{IP}(0)$ are lengths of arcs from $E(0)$ and $P(0)$ to $I$, respectively,
while $\overarc{IE}(0)/\mu v_P$ and $\overarc{IP}(0)/v_P$ are the corresponding times for $E$ and $P$ 
to reach $I$ using their maximal speeds.

\begin{figure}
\begin{center}
\includegraphics[width=0.4\textwidth]{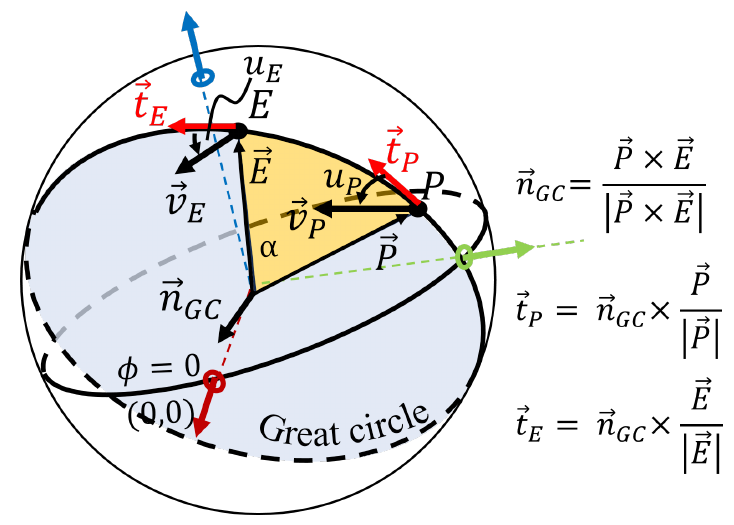}
\end{center}
\caption{Relative position between $P$ and $E$: The unit vector $\vec{n}_{GC}$ defines the orientation or the great circle passing through $P$ and $E$. The relative $P$-$E$ position is given by the angular distance $\alpha$ measured along the great circle in the direction defined by $\vec{n}_{GC}$ and the right hand rule. The $P$ control variable is the angle $u_P$, $E$ control variables are the angle $u_E$ and speed $v_E=|\vec{v}_E|$. The angles $u_E$ and $u_P$ are measured with respect to the great circle, i.e., with respect to the tangent vector $\vec{t}_E$ and $\vec{t}_P$.}\label{Fig2}
\end{figure}
\section{Pursuit-Evasion equilibrium strategies}
\label{sec:strategies}
Figure~\ref{Fig2} illustrates that the center of sphere, points $E$ and 
$P$ define a plane $\mathcal{G}_C$ with the unit vector $\vec{n}_{GC}$ 
defined as 
\begin{equation}
 \vec{n}_{GC} \sin \alpha = \tfrac{1}{R^2}\vec{P} \times \vec{E}.
\label{eq:nGCalpha}
\end{equation}
The \emph{angular distance}, $\alpha$, is the angle between $\vec{P}$ and $\vec{E}$ 
and measured in the direction defined by $\vec{n}_{GC}$ 
and the right hand rule.

The intersection of the plane $\mathcal{G}_C$ and the sphere's surface defines a great circle which passes through $E$ and $P$. \emph{The $P-E$ relative angular position $\alpha$ is described by the orientation of the great circle given by the vector $\vec{n}_{GC}$ and non-negative angular distance $\alpha \ge 0$ along the great circle}. Note that for co-linear vectors $\vec{P}$ and $\vec{E}$, the vector $\vec{n}_{GC}$ is not well defined, meaning that the direction in which we measure $\alpha$ is not well defined. However, we can still define $\alpha \in [0,\pi]$ to include these cases as
\begin{equation*}
\alpha = \begin{cases}
    \arccos \tfrac{1}{R^2} (\vec{P} \cdot \vec{E}),& \vec{P} \times \vec{E} \ne 0 \\
    \pi,& \vec{P}=-\vec{E} \\
    0,& \vec{P}=\vec{E},
\end{cases}
\end{equation*}
where $\cdot$ is a scalar product between $\vec{P}$ and $\vec{E}$.
Given \eqref{eq:veE}, \eqref{eq:vpP} and \eqref{eq:nGCalpha}, the velocities of the evader and pursuer are
\begin{align}
   \vec{v}_E&=v_E \vec{e}= v_E \overbrace{(\cos u_E \vec{t}_E + \sin u_E \vec{n}_{GC})}^{\vec{e}} \label{eq:evader}\\
   \vec{v}_P&=v_P\vec{p}= v_P \underbrace{(\cos u_P \vec{t}_P + \sin u_P \vec{n}_{GC})}_{\vec{p}},
    \label{eq:pursuer}
\end{align}
using the unit vector $\vec{n}_{GC}$ and  vectors $\vec{t}_E$, $\vec{t}_P$ that are tangent to the great circle at points $E$ and $P$. In these expressions, $u_E$ and $u_P$ are angles between the corresponding velocity vectors and tangent to the great circle at $E$ and $P$, respectively. Consequently, instead of the vector control variable $\vec{e}$ for $E$, we can use the scalar variable $u_E \in [0, 2\pi)$, and instead of the vector control variable $\vec{p}$ for $P$, we can use the scalar variable $u_P \in [0, 2\pi)$. 

Using $u_E$ and $u_P$, we show (see Appendix) for the kinematics of $E$ and $P$ from \eqref{eq:evader} and \eqref{eq:pursuer}, respectively, that the relative angular distance $\alpha > 0$ obeys
\begin{equation}
   \dot{\alpha} = \tfrac{v_E}{R} \cos u_E - \tfrac{v_P}{R} \cos u_P, \label{eq:relativPos}
\end{equation}
for all $\alpha \in (0, \pi)$ when the great circle (GC) is well defined. The special case $\alpha=0$ corresponds to the capture. The other special case $\alpha=\pi$ will be considered as a part of the game analysis in Lemma~\ref{lem:dispersal-control}.

Let us assume that the initial relative position between $P$ and $E$ is $\alpha(0) \in(0,\pi]$ and define a cost function 
\begin{equation}
J =\int_{0}^\tau 1 \cdot \mathrm{d}t = \tau \label{eq:Cost},
\end{equation}
where $\tau$ is the capture time at which $\alpha(\tau)=0$. 
The goal of $P$ is to minimize $J$, i.e., time to capture and 
the goal of $E$ is to maximize it.  Theorem~\ref{thm:1v1-Value} defines equilibrium strategies for $P$ and $E$ and their relative angular distance $\alpha\in (0,\pi)$. 
This is followed by Lemma~\ref{lem:dispersal-control} about the case $\alpha=\pi$ and all 
the results are summarized in Theorem~\ref{thm:full-game-solution} providing equilibrium strategies for $P$-$E$ relative positions $\alpha \in (0,\pi]$.

\begin{theorem}
\label{thm:1v1-Value}

For the P-E game on a sphere with agent 
kinematics \eqref{eq:evader}-\eqref{eq:pursuer}, time to 
capture cost \eqref{eq:Cost} and 
\begin{itemize}
\item[($a$)] The relative angular position between $P$ 
and $E$ is $\alpha~\in~(0,\pi)$ and obeys \eqref{eq:relativPos} 
\item[($b$)] $P$ has constant speed $v_P$, the control variable is heading $u_P\in [0, 2\pi)$ and the goal to minimize the time to capture
\item[($c$)] The control variables for $E$ are speed $v_E \in [0, \mu v_P]$ 
and heading $u_E \in [0, 2\pi)$, and $E$ has the goal to maximize the time 
to capture,
\end{itemize}
the saddle point equilibrium strategies for all 
$\alpha \in (0,\pi)$  are defined by 
\begin{equation}
u_P = 0\ \ \ {\rm and}\ \ \  (v_E,u_E) = 
  (\mu v_P, 0). \label{eq:optStrat}
\end{equation}
Moreover, if both players are at the relative 
angular position $\alpha \in (0,\pi)$ and follow their equilibrium strategies, the 
time to capture is 
\begin{equation}
V(\alpha)= R\alpha ((1-\mu)v_P)^{-1}. 
\label{eq:value}
\end{equation}
\end{theorem}

\begin{proof}
The equilibrium strategies satisfy the Hamilton-Jacobi-Isaacs (HJI) equation
\begin{equation}
\min_{u_P} \max_{(v_E, u_E)}\left\{\left(\tfrac{v_E}{R} \cos u_E - 
\tfrac{v_P}{R} \cos u_P\right) \frac{\partial V}{\partial \alpha} + 
1 \right\}=0, \label{eq:HJI1}
\end{equation}
where $V(\alpha)$ is the Value of the game (see Remark~\ref{rem:Value}),
the minimization is over the variable $u_P$ defining the heading of $P$ and the maximization is over the speed 
$v_E$ and variable $u_E$ defining the heading of $E$. When 
the $\min$ and $\max$ operators are applied, we obtain
\begin{equation}
 \max_{(v_E, u_E)}\left\{\tfrac{v_E}{R} \cos u_E \tfrac{\partial V}{\partial \alpha} \right\} +\min_{u_P}\left\{- \tfrac{v_P}{R} \cos u_P \tfrac{\partial V}{\partial \alpha}\right\} + 1 =0. \label{eq:HJI2}
\end{equation}

The fact that $\max$ and $\min$ operators are separable reinforces that the strategies in \eqref{eq:HJI1} are saddle point strategies.
The game is defined over variable $\alpha$ and since $\alpha=0$ corresponds to the capture, we know that $V(0)=0$. Also, we know that for an incrementally larger angular distance $\alpha \in (0,\pi)$, the game lasts incrementally longer, i.e., $\frac{\partial V}{\partial \alpha} > 0$. Therefore, we can conclude
\begin{itemize}
\item  $\min$  in the second term of \eqref{eq:HJI2} is achieved for $u_P=0$, which is
the equilibrium $P$ strategy  for $\alpha \in (0,\pi)$
\item $\max$ in the first term of \eqref{eq:HJI2} is achieved for $v_E = \mu v_P$, $u_E=0$, which is the equilibrium $E$ strategy for $\alpha \in (0,\pi)$. 
\end{itemize}
\emph{These facts prove the equilibrium actions in} 
\eqref{eq:optStrat}. 

When we substitute the equilibrium strategies in \eqref{eq:HJI2}, we obtain 
\begin{equation*}
\tfrac{\mu v_P}{R} \tfrac{\partial V}{\partial \alpha} - \tfrac{v_P}{R} \tfrac{\partial V}{\partial \alpha} + 1 =0. %\label{eq:HJI2}
\end{equation*}
Re-arranging, we obtain
\begin{equation*}
\tfrac{\partial V}{\partial \alpha}= \tfrac{R}{(1-\mu)v_P}.
\end{equation*}
Finally, since $V(0)=0$, we can conclude that 
\begin{equation*}
V(\alpha)= \tfrac{R\alpha}{(1-\mu)v_P},
\end{equation*}
which is \eqref{eq:value} and that concludes our proof.%\qed \\
\end{proof}

\begin{remark}
% \emph{Remark~1.}
\label{rem:Value}
The cost function resulting from both agents following their equilibrium strategies is also called the Value of the game, hence the symbol $V$ in \eqref{eq:value}. 
\end{remark}

In Lemma~\ref{lem:dispersal-control}, we discuss the special case when the game starts 
from the initial position $\alpha=\pi$. In this case, the great 
circle is not well defined and there are infinitely many great 
circles that go through $P$ and $E$, yet the angular distance 
$\alpha$ measured along any of them is $\alpha=\pi$. \\

\begin{remark}
The essential problem with $\alpha=\pi$ is 
that without a well-defined great circle, there are no reference 
directions for measuring angles $u_E$ and $u_P$.
\end{remark}

\begin{lemma}
    \label{lem:dispersal-control}
For all conditions of Theorem~\ref{thm:1v1-Value}, except 
($a$) which is replaced with ($a'$): the relative position 
$\alpha=\pi$, the equilibrium actions for $P$ and $E$ are
\begin{align}
u_P &= {\mathrm{any}\ \mathrm{direction}} \\
(v_E,u_E) &=(0,u_E), \quad u_E \in [0, 2\pi), \label{eq:stratPI}
\end{align}
where $u_E$ is measured with respect to a great circle 
aligned with the direction of motion selected by $P$. 
The equilibrium action for $E$ corresponds to 
\begin{equation*}
   \min_{v_E } \max_{u_E} \mathcal{L}(v_E,u_E)= 
   \max_{u_E } \min_{v_E} \mathcal{L}(v_E,u_E)=\tfrac{\mu}{1-\mu}, 
   \label{eq:minmax}
\end{equation*}
\begin{equation*}
\mathcal{L}(v_E,u_E)=\frac{\mu-\frac{v_E}{v_P} \cos u_E}{1-\mu}, 
v_E \in [0, \mu v_P], u_E \in [0, 2\pi). 
\end{equation*}
\end{lemma}

\begin{proof}
% \emph{Proof}.
Due to infinitely many great circles that we can consider to 
measure $\alpha$, we are free to take the one which is aligned with 
the selected direction of motion by $P$. Because of that, whatever 
direction is selected by $P$, its $u_P=0$ since the adopted great 
circle is aligned with its motion, but $u_E$ has to be also measured 
with respect to the same great circle. The relative angle on the adopted 
great circle for $\alpha=\pi$ satisfies
\begin{equation}
\dot{\alpha} = \tfrac{v_E}{R} \cos u_E - \tfrac{v_P}{R} < 0,  
\label{eq:dotalpha_pi}
\end{equation}
which follows from $v_E \le \mu v_P$ and $0 < \mu < 1 \Rightarrow v_E < v_P$. 
In other words, on the adopted great circle where $u_P=0$ and for $\alpha(0)=\pi$, 
any action of $E$ results in 
\begin{equation}
(v_E, u_E) \in \{[0, \mu v_P] \times [0,2\pi)\} \Rightarrow \alpha(0^+) < \pi, 
\label{eq:hplus}
\end{equation}
and we can conclude that at $t=0^+$, the great circle is well defined. 
Consequently, for any $t>0^+$, both $P$ and $E$ can follow their equilibrium 
strategies from Theorem~\ref{thm:1v1-Value} since $\alpha(t) < \pi$ for all $t>0^+$.

However, for $t=0$ and the adopted great circle ($u_P=0$) in \eqref{eq:dotalpha_pi}, 
we can only state that $u_E \in [0,2\pi)$ which is measured with respect to the 
adopted great circle. It is because the adopted great circle, i.e., the direction of motion 
selected by $P$ is unknown to $E$. Therefore, whatever heading direction is selected 
by $E$, the best we can say is that $u_E \in [0,2\pi)$. To find equilibrium actions at 
$t=0$, we will use the rate-of-loss analysis \cite{TAC2023} which starts 
with the following reasoning. 

Let us assume that at $t=0$, $\alpha(0)=\pi$ and both $P$ and $E$ commit 
their control variable selections that hold for a very short (hold time) 
time $h \approx 0$, $h>0$. As was previously discussed, due to the adopted great circle at $t=0$, we always 
have $u_P=0$; therefore after $h$, the value $\alpha(h)$ is
\begin{equation*}
\alpha(h) = \pi + \frac{v_E h \cos u_E}{R}- \frac{v_P h}{R} < \pi - \frac{v_P-v_E}{R}h < \pi.
\end{equation*}
Since $\alpha(h)<\pi$, Theorem~\ref{thm:1v1-Value} provides strategies for both $P$ and $E$, and the resulting time to capture is $h+V(\alpha(h))$, which depends on where $\alpha$
landed after the initial hold time $h$. In general, $h+V(\alpha(h))$ can be different 
from $V(\alpha(0))$. To measure that discrepancy as $h \rightarrow 0$, we 
use the instantaneous rate-of-loss   
\begin{align}
\mathcal{L} & = - \lim_{h \rightarrow 0} \frac{V(\alpha(h))-V(\alpha(0))+h}{h} 
\label{eq:rateLoss} \\
&=-\lim_{h \rightarrow 0}\tfrac{1}{h}\left(\frac{R(\pi+\frac{v_E h \cos u_E-v_Ph}{R})}{(1-\mu)v_P}
- \frac{R\pi}{(1-\mu)v_P} + h \right) \nonumber \\
&=-\frac{v_E  \cos u_E-v_P}{(1-\mu)v_P}-1 = \tfrac{1}{1-\mu}\left(\mu-\frac{v_E}{v_P} \cos u_E \right) \ge 0. \ \label{eq:rateofloss2}
\end{align}
This is the rate-of-loss from the prospective of $E$. Its loss is reflected in the shortening 
of its time to capture, hence we  included the ``--" sign in front of the limit of \eqref{eq:rateLoss} to be able to discuss the rate-of-loss in terms of a non-negative value.

In expression \eqref{eq:rateofloss2}, the rate of loss $\mathcal{L}=\mathcal{L}(v_E, u_E)$ and for any $v_E \in [0, \mu v_P]$ the maximum of loss $\mathcal{L}(v_E,u_E)$ corresponds to $u_E=-\pi$, i.e., 
\begin{equation*}
\max_{u_E} \mathcal{L}(v_E, u_E)=\tfrac{1}{1-\mu}\left(\mu + 
\frac{v_E}{v_P}\right).
\end{equation*}
We focus on this max of the rate-of-loss as the worst case
scenario since $u_E$ does not know the selected direction of movement by $P$, i..e, the orientation of the adopted 
great circle. This max increases with $v_E$, so the only 
way that $E$ can minimize it is to select $v_E=0$ and 
we can conclude that 
\begin{equation}
\min_{v_E} \max_{u_E} \mathcal{L}(v_E, u_E)=\tfrac{\mu}{1-\mu} \label{eq:lossrate}
\end{equation}
corresponds to action $(0, u_E)$, $u_E=[0, 2\pi)$, which proves 
the property of the equilibrium action for $E$ in \eqref{eq:stratPI}.

Alternately, rather than first minimizing over $u_E$, if we initially apply minimization over $v_E$ from 
\eqref{eq:rateofloss2}, we obtain 
\begin{equation*}
\min_{v_E} \mathcal{L}(v_E, u_E)=
\begin{cases}
    l_1=\frac{\mu-\frac{\mu v_P}{v_P} \cos u_E}{1-\mu}, &\cos u_E \ge 0 \\ 
\l_2=\frac{\mu-\frac{0}{v_P} \cos u_E}{1-\mu}, & \cos u_E < 0. 
\end{cases}
\end{equation*}
Then, when we apply max over $u_E \in [0,2\pi)$, we see that $l_1$ reaches 
the max value $l_1=\tfrac{\mu}{1-\mu}$, but only if $u_E=\pi/2$. Since $E$ 
does not know the orientation of the adopted great circle, it cannot 
use that value. On the other hand, the max value $l_2=\tfrac{\mu}{1-\mu}$
is reached not only for $\cos u_E <0$, but for any angle $u_E \in [0, 2\pi)$. 
From this and \eqref{eq:lossrate}, we can finally conclude that the equilibrium 
action for $E$ which is $(0, u_E)$, $u_E \in [0, 2\pi)$ corresponds to the saddle point
\begin{align}
\min_{v_E} \max_{u_E} \mathcal{L}(v_E,u_E)=\max_{u_E} \min_{v_E} \mathcal{L}(v_E,u_E) =\tfrac{\mu}{1-\mu},
\end{align}
which proves \eqref{eq:minmax} and with this, we conclude our proof.
\end{proof}

We finish this section with Theorem~\ref{thm:full-game-solution} which summarizes results 
from Theorem~\ref{thm:1v1-Value} and Lemma~\ref{lem:dispersal-control}.
\begin{theorem}
\label{thm:full-game-solution}
% \emph{Theorem~\ref{thm:full-game-solution}}.
For all conditions of Theorem~\ref{thm:1v1-Value} except ($a$) 
which is replaced with $(a'')$: all possible relative positions
$\alpha~\in~(0, \pi]$, the equilibrium actions for $P$ and $E$ are
\begin{align}
u_P &= \begin{cases}
  0, & \alpha \in (0,\pi) \\
  {\mathrm{any}\; \mathrm{direction}},& \alpha = \pi  
\end{cases}
\label{eq:pactThm2}\\
(v_E,u_E) &= \begin{cases}
  (\mu v_P, 0),&  \alpha \in (0,\pi) \\
  (0,{\mathrm{any}\; \mathrm{direction}}),& \alpha = \pi.  
\end{cases}
\label{eq:eactThm2}
\end{align}

For any initial condition $\alpha \in (0,\pi]$, and 
$P$ and $E$ following their equilibrium strategies, the time to capture is 
\begin{equation}
\tau \le V(\alpha) = \tfrac{R \alpha}{(1-\mu)v_P}, \quad \mathrm{for} \; \alpha \ne \pi,
\label{eq:valThm2}
\end{equation}
%in which ``=" holds for $\alpha \ne \pi$. 
and for $\alpha=\pi$ the difference 
\begin{equation}
V(\pi)-\tau = \varepsilon >0,  
\label{eq:varepsilon}
\end{equation}
where $\varepsilon$ is an infinitesimally small value, $\varepsilon \rightarrow 0$.
\end{theorem}

\begin{proof}
% \emph{Proof}.
The equilibrium actions in \eqref{eq:pactThm2} and \eqref{eq:eactThm2} 
follow directly from Theorem~\ref{thm:1v1-Value} and Lemma~\ref{lem:dispersal-control} covering the cases $\alpha \in (0, \pi)$ 
and $\alpha = \pi$, respectively.  In \eqref{eq:stratPI}, Lemma~\ref{lem:dispersal-control} states that 
$u_E \in [0, 2\pi)$ and here we restate it as $u_E=$ any direction since the range for $u_E$ covers all possible angles. 
Also note that the angle $u_E$ is measured with respect to the great circle aligned with the direction of movement selected by $P$, which can also take any direction. 
This justifies that instead of referring to a range, Theorem~\ref{thm:full-game-solution} states that $u_E=$ any direction. 
With this, we conclude the discussion of \eqref{eq:pactThm2} and \eqref{eq:eactThm2} in this proof.

With regard to \eqref{eq:valThm2}, in Theorem~\ref{thm:1v1-Value} we already proved that for 
$\alpha \in (0,\pi)$ the equality in \eqref{eq:valThm2} holds. Therefore, we 
only need to prove that for $\alpha=\pi$ the time to capture $\tau$ satisfies  
\begin{equation}
 \tau \le V(\pi) = \tfrac{R \pi}{(1-\mu)v_P}. 
 \label{eq:prove}
\end{equation}
For $\alpha=\pi$ in \eqref{eq:rateLoss}--\eqref{eq:rateofloss2} of the Lemma~\ref{lem:dispersal-control} proof, 
we use the rate-of-loss 
\begin{align}
\mathcal{L}(v_E, u_E)&=- \lim_{h \rightarrow 0} \frac{V(\alpha(h))-V(\alpha(0))+h}{h} \label{eq:rateoflossa} \\
&= \tfrac{1}{1-\mu}\left( \mu-\tfrac{v_E}{v_P} \cos u_E \right), \label{eq:rateofloss3b}
\end{align}
where $h>0$ is an infinitesimal hold time, $h \rightarrow 0$. Let us rewrite 
\eqref{eq:rateoflossa} for $\alpha(0)=\pi$ as 
\begin{equation}
V(\pi)-(V(\alpha(h))+h)  = h \mathcal{L}(0, u_E) , h \rightarrow 0, 
\label{eq:vpitau}
\end{equation}
and recognize that for the equilibrium action of $P$, $\alpha(h) < \pi$ which follows 
from (\ref{eq:hplus}). Furthermore, if both $P$ and $E$ follow their equilibrium 
strategies for all $t>h$, $V(\alpha(h))$ is the exact time until the capture. Therefore, 
the time to capture from $t=0$ is $\tau=V(\alpha(h))+h$ which we substitute 
in \eqref{eq:vpitau} together with $\mathcal{L}(0, u_E)=\tfrac{\mu}{1-\mu}$ to 
obtain 
\begin{equation}
  \tau  = V(\pi)-h \tfrac{\mu}{1-\mu} , \; h \rightarrow 0 \Rightarrow \tau \le V(\pi),
\end{equation}
which proves \eqref{eq:prove} and consequently \eqref{eq:valThm2}. The same 
expression yields
\begin{equation}
     V(\pi)-\tau = \underbrace{h \tfrac{\mu}{1-\mu}}_{\varepsilon}, \; h \rightarrow 0, 
\end{equation}
where the $\varepsilon$ expression on the right is an infinitesimally small value, 
which proves \eqref{eq:varepsilon} and concludes our proof. 
% \qed
\end{proof}

\section{Apollonius Domain and Intercept Point}
\label{sec:Apollonius}

\begin{figure}
\begin{center}
\includegraphics[width=0.4\textwidth]{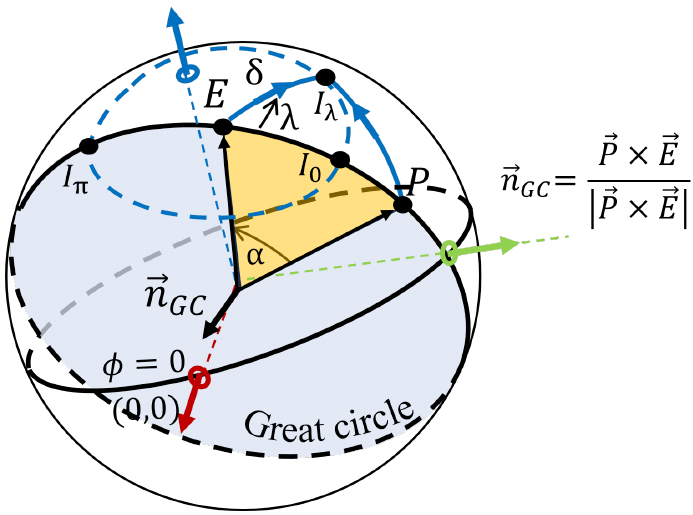}
\end{center}
\caption{The Apollonius domain $\mathcal{A}$ is depicted by a set of points enclosed in the blue dashed curve connecting points $I_0$, $I_\lambda$ and $I_{\pi}$. 
The blue dashed curve is the boundary of Apollonius domain, $\partial \mathcal{A}$, and it belongs to the domain as well. 
The Apollonius domain includes $E$ and its boundary is a set of intercept points $I_{\lambda}$ satisfying $|\overarc{PI}_{\lambda}|/\delta(\lambda)=\mu$, where $\delta(\lambda)$ and $|\overarc{PI}_{\lambda}|$ are distances from $E$ and $P$ to $I_{\lambda}$, respectively. 
There is one-to-one correspondence between the angle $\lambda$ and the distance $\delta(\lambda)$.}\label{Fig3}
\end{figure}

Figure~\ref{Fig3} depicts $P$ and $E$ at a relative angular position $\alpha$. 
Let us assume that at initial time $t=0$, $E$ picks a direction of motion defined by an angle $\lambda$ and starts moving with its maximal speed $v_E=\mu v_P$. 
The angle $\lambda$ is measured with respect to the great circle defined by $P$ and $E$, and it is known by $P$ at $t=0$ as well. 
Using this information, $P$ can take the \emph{shortest path} to intercept $E$ at the point $I_{\lambda}$. 
The arc length from $E$ to $I_{\lambda}$ is $\delta =|\overarc{EI}_{\lambda}|=\mu v_P \tau_{\lambda}$ while the arc length from $P$ to $I_{\lambda}$ is $|\overarc{PI}_\lambda|=v_P \tau_{\lambda}$, where $\tau_{\lambda}$ is the time of the intercept. 
Note that all points on the arc $\overarc{EI}_\lambda$ are points that $E$ can reach sooner than $P$; similarly, $\overarc{PI}_{\lambda}$ are points that $P$ can reach sooner than $E$. 
We can sweep the angle $\lambda \in [0, 2\pi)$ over its full range, and it will produce a closed line of $I_{\lambda}$ points. 
Points inside (outside) of the line are all points that $E$ ($P$) can reach  sooner. 
The points inside the line form the Apollonius domain $\mathcal{A}$ and the boundary of the domain $\partial \mathcal{A}$ is a line which $P$ and $E$ can reach simultaneously, therefore, it is a set of potential intercept points from the initial configuration.

It is shown in \cite{Kovshov2000} that 
\begin{align}
\cos \left( \tfrac{\delta}{R\mu}\right)  = \cos \tfrac{\delta}{R} \cos \alpha +  \sin \alpha \cos \lambda \sin \tfrac{\delta}{R}, 
    \label{eq:lambdadeltarel}
\end{align}
and that this implicit relation describes a one-to-one $\delta(\lambda)$ mapping between 
$\lambda$ and $\delta$ which is dependent on $\alpha$. It is shown in \cite{Kovshov2000} that 
there is a critical value $\alpha_c =\pi(1-\mu)$ such that 
for $\lambda \in [0, \pi]$ and $\alpha \le \alpha_c$ 
\begin{align}  
\delta(0)=
\tfrac{\alpha R\mu}{1+\mu} \le \delta(\lambda) \le \tfrac{ \alpha R\mu}{1-\mu} = \delta(\pi), \label{eq:delta_lt_alphac}
\end{align}
and for $\lambda \in [0, \pi]$ and $\alpha > \alpha_c$ 
\begin{align}    
\delta(0)=\tfrac{\alpha R\mu}{1+\mu} \le \delta(\lambda) \le \tfrac{ R \mu(2\pi-\alpha)}{1-\mu} =\delta(\pi). \label{eq:delta_gt_alphac}
\end{align}
In both cases, the lower bound corresponds to $\delta(0)=\overarc{EI}_0$ which is the arc length between $E$ and $I_0$, and the upper bound corresponds to $\delta(\pi)=\overarc{EI}_{\pi}$ which is the arc length between $E$ and $I_{\pi}$. 
Note that \eqref{eq:lambdadeltarel} can be rewritten as
\begin{equation}
\cos \lambda = \left(\cos \left( \tfrac{\delta}{\mu R}\right) - \cos \tfrac{\delta}{R} \cos \alpha \right) \left(\sin \alpha  \sin \tfrac{\delta}{R} \right)^{-1},
\end{equation}
and due to this  
\begin{equation}
\delta(\lambda)=\delta(-\lambda),\ {\rm for} \ \lambda \in [-\pi,0]. 
\label{eq:lambdasymmetry}
\end{equation}
Expressions \eqref{eq:lambdadeltarel}-\eqref{eq:lambdasymmetry} are sufficient to define 
the boundary of Apollonius domain, $\partial \mathcal{A}$, for any  
initial configuration between $P$ and $E$.

\begin{theorem}
\label{thm:intercept-point-on-apollonius-domain}
For all conditions of \emph{Theorem~\ref{thm:full-game-solution}} and if both $P$ and 
$E$ use their equilibrium strategies, the intercept point $I_{P-E}$ belongs 
to the boundary of Apollonius domain, $\partial \mathcal{A}$, only if the initial 
relative distance $\alpha(0)$ between $P$-$E$ satisfies 
\begin{equation}
 \alpha(0) \le \alpha_c= \pi(1-\mu). \label{eq:condThm3}
\end{equation}
\end{theorem}

\begin{proof}
Since $\mu < 1$, the condition \eqref{eq:condThm3} implies 
that $\alpha(0)<\pi$. Therefore, the equilibrium actions for both $E$ and 
$P$ are $u_E=0$ and $u_P=0$. With $u_E=0$, $E$ moves in  
the direction $\lambda=\pi$ and the intercept point is on the great circle 
defined by $P$ and $E$. The time to capture is (see Theorem~\ref{thm:full-game-solution})
\begin{equation}
\tau = V(\alpha(0)) = \tfrac{R \alpha(0)}{(1-\mu)v_P}. 
\end{equation}
Therefore, the arc length traveled by $E$ is 
\begin{equation}
\overarc{EI}_{P-E}=\tfrac{R \mu v_P \alpha(0)}{(1-\mu)v_P}= \tfrac{R \mu  \alpha(0)}{(1-\mu)}, \label{eq:arcEIPE}
\end{equation}
which is a value independent of $\alpha_c$. Now, for $\alpha(0)<\alpha_c$ from (\ref{eq:delta_lt_alphac}), we have 
\begin{equation}
\overarc{EI}_{\pi}=\delta(\pi)=\tfrac{R \mu  \alpha(0)}{(1-\mu)} =\overarc{EI}_{P-E}, 
\end{equation}
and we can conclude that $I_{P-E}=I_{\pi}$. However, for $\alpha(0)>\alpha_c$ from (\ref{eq:delta_gt_alphac}), we have 
\begin{equation}
\overarc{EI}_{\pi}=\delta(\pi)=\tfrac{R \mu(2\pi-\alpha(0))}{(1+\mu)}.  
\end{equation}
If we assume that $I_{P-E}=I_{\pi}$, then from (\ref{eq:arcEIPE}) 
\begin{equation}
\overarc{EI}_{\pi}=\tfrac{R(2\pi-\alpha(0))}{(1+\mu)}= \tfrac{R\alpha(0)}{(1-\mu)}=\overarc{EI}_{P-E},
\end{equation}
which yields
\begin{equation}
\alpha(0)=\pi(1-\mu)=\alpha_c,
\end{equation}
and contradicts $\alpha(0) > \alpha_c$. Therefore, $I_{P-E} \ne I_{\pi}$, which concludes the 
proof of the theorem. %\qed \\
\end{proof}

\section{Results}
\label{sec:Results}

\begin{figure}[t]
    \centering
    \subfloat[
        $v_E=0.35$ %, {\color{red}$α=40 < α_C$}, {\color{green}$α = α_C = 72$}, {\color{blue}$α = 80 > α_C$}
    ]{\includegraphics[trim={1cm 1cm 1cm 1cm},clip,width=0.49\linewidth]{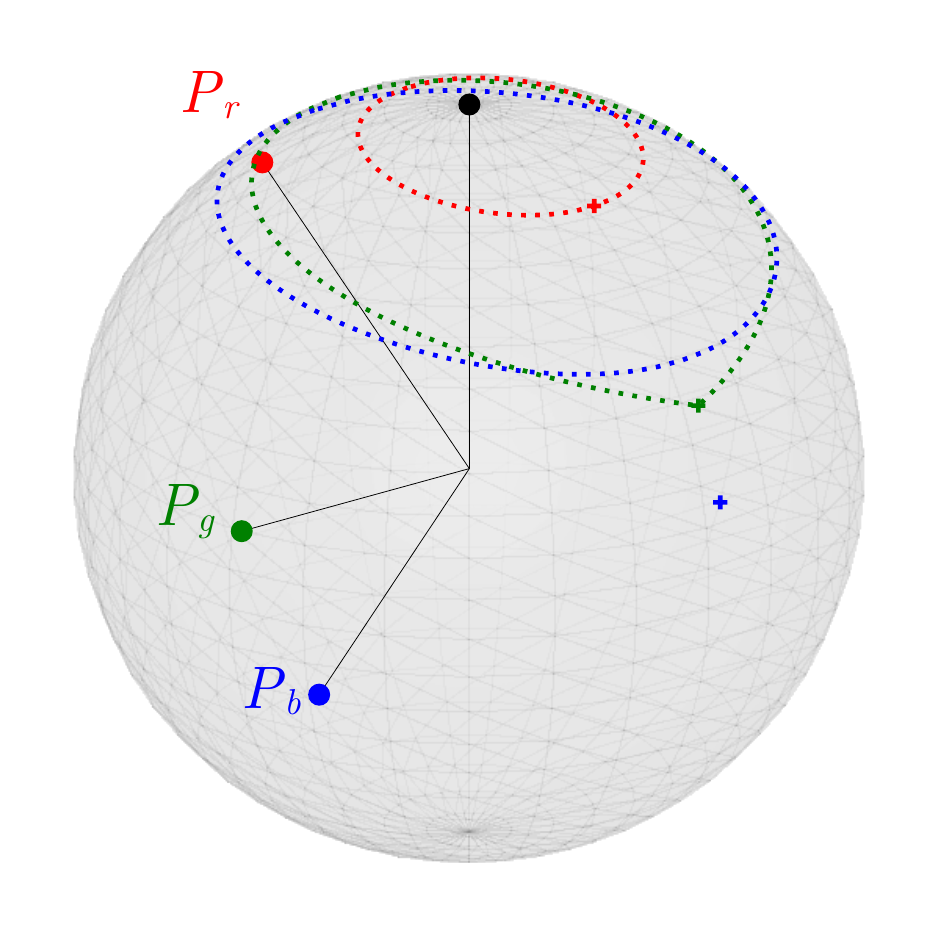}
    \label{fig:Res1}}
    \subfloat[
        $v_E=0.6$%, {\color{red}$α=80 < α_C$}, {\color{green}$α = α_C=117$}, {\color{blue}$α = 140 > α_C$}
    ]{\includegraphics[trim={1cm 1cm 1cm 1cm},clip,width=0.49\linewidth]{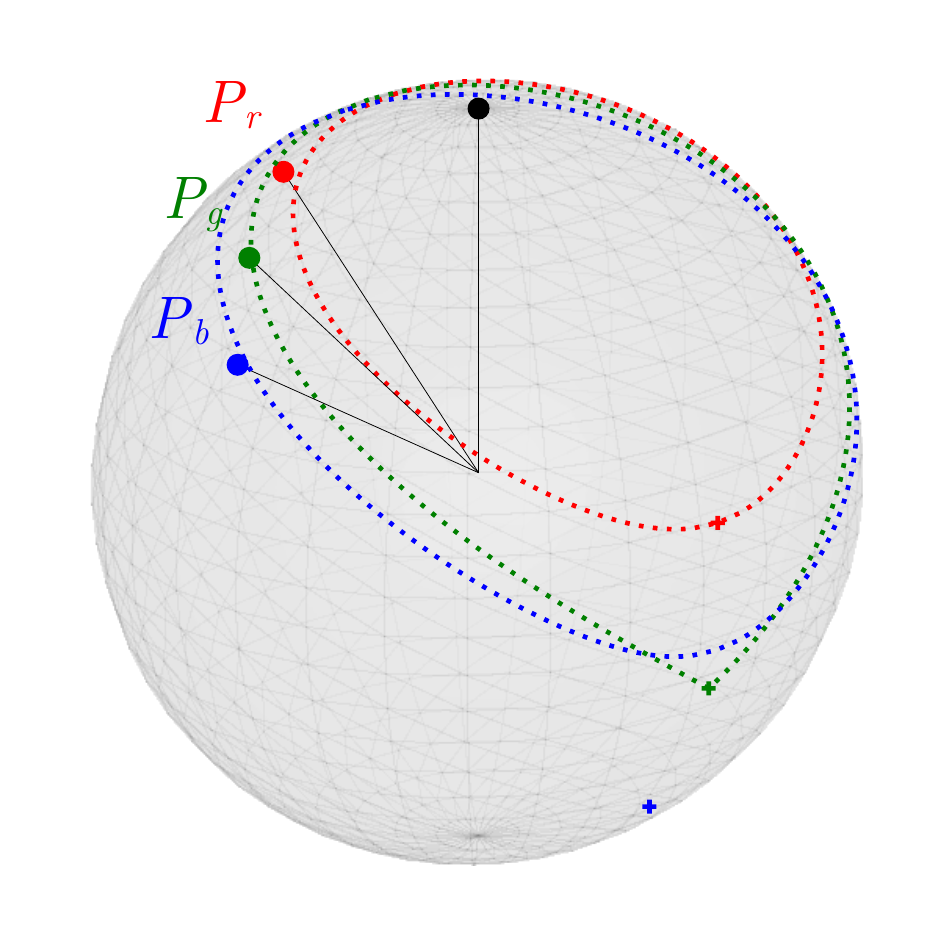}
    \label{fig:Res2}}
    \caption{
        Apollonius domain for different $v_E$ (with $v_P = 1$) and $α$:
        The colors (also indicated by the subscripts) correspond to {\color{red}$\alpha_r < \alpha_C$}, {\color{green!50!black}$\alpha_g = \alpha_C$}, and {\color{blue}$\alpha_b > \alpha_C$}.
        The dashed lines indicate the Apollonius domains and the small $+$ indicates the equilibrium intercept point, $I_{P-E}$.
    }
    \label{fig:AD-various-parameters}
\end{figure}

Figure~\ref{fig:AD-various-parameters} illustrates the shape of the Apollonius domain for two different evader speeds and for initial geodesic distances below, at, and above the critical $α_C$ value from Theorem~\ref{thm:intercept-point-on-apollonius-domain}.
In addition, the equilibrium intercept point for the game of $\min$-$\max$ capture time is drawn.
As expected, the intercept point lies on the associated Apollonius domain when $α \leq α_C$ and lies outside when $α > α_C$. Towards the end of this 
section, we present applications of our results to two  differential 
game problems on the sphere with more than one $P$.

\subsection{Two Pursuers}
\label{sec:Two_Pursuers}

In the spirit of Isaacs~\cite[Example 6.8.3]{Isaacs1965}, we now extend the results for the one-pursuer scenario to the case of two cooperative pursuers by leveraging the geometry of the Apollonius domain.
First let us define the evader's dominance region as $\mathcal{E} = \mathcal{A}_1 \cap \mathcal{A}_2$, where $\mathcal{A}_1$ and $\mathcal{A}_2$ are the Apollonius domains associated with pursuers $P_1$ and $P_2$, respectively. 
The set $\mathcal{E}$ represents the set of points that the evader can reach before either pursuer.
With the evader at the north pole, let $α_1$ and $α_2$ denote the latitude offset for the pursuers, and let $λ_o$ be the longitudinal offset between the pursuers.
The agents' speeds are such that $μ_1 = \tfrac{v_E}{v_{P_1}}, μ_2 = \tfrac{v_E}{v_{P_2}} < 1$.
Let the point $I_{P_1-E}$ be the equilibrium intercept point for the one-on-one game between $P_1$ and $E$ (and let $I_{P_2-E}$ be defined similarly for $P_2$).

\begin{proposition}
    \label{prop:2v1}
    In the two-pursuer scenario with $α_1 < π (1 - μ_1)$ and $α_2 < π (1 - μ_2)$ (i.e., the pursuers are sufficiently close to the evader), the intercept point associated with the solution for the game of $\min$-$\max$ capture time is
   \begin{equation}
    \label{eq:two-v-one-solution}
    I^* =
    \begin{cases}
        I_{P_1-E}, & I_{P_1-E} \in \mathcal{A}_2  \\
        I_{P_2-E}, & I_{P_2-E} \in \mathcal{A}_1  \\
        \max \left\{ \overarc{EI} \mid I \in \partial\mathcal{A}_1 \cap \partial\mathcal{A}_2\right\}, & \mathrm{otherwise.}
    \end{cases}
   \end{equation}
\end{proposition}

\begin{figure}[htbp!]
    \centering
    \includegraphics[trim={1cm 1cm 1cm 1cm},clip, width=0.64\linewidth]{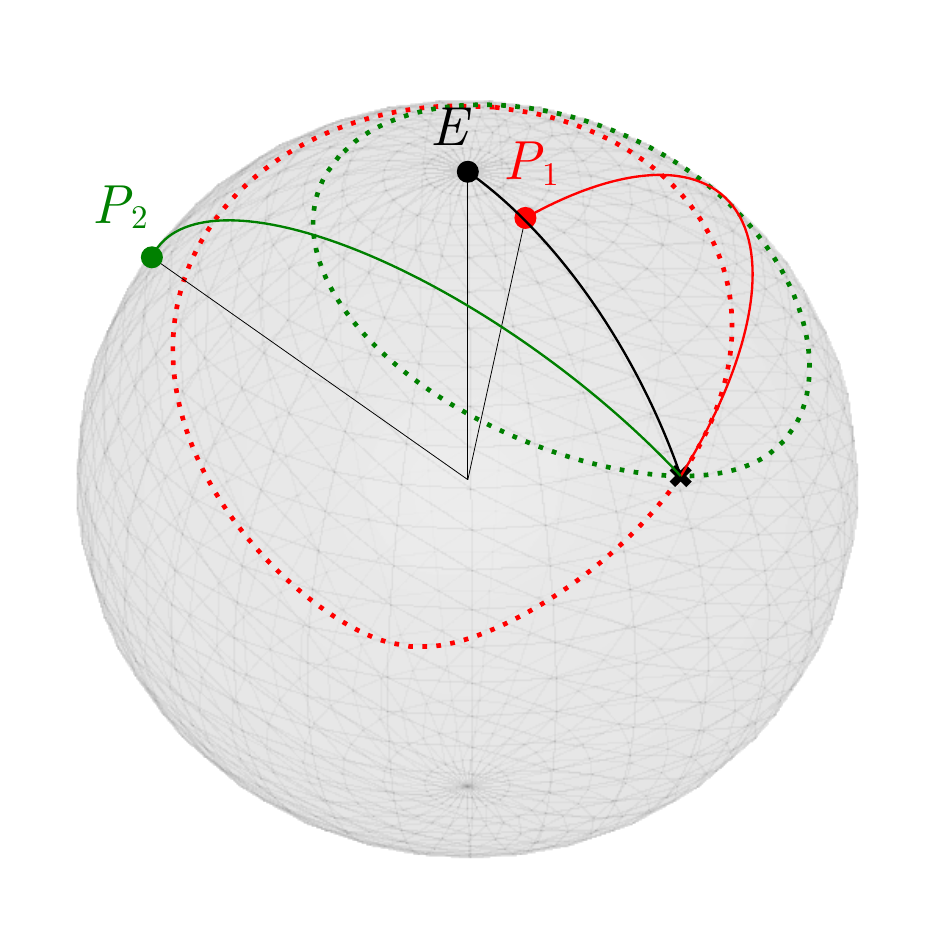}
    \caption{
        Proposed solution for the two-pursuer, one-evader $\min$-$\max$ capture time game with $μ_1 = μ_2 = 0.5$, $α_1 = 0.9 π (1 - μ) \approx 1.41$, $α_2 = 0.8 α_1 \approx 1.13$, and with the pursuers' longitudinal offset as $λ_o = 0.4 π \approx 1.26$.
        The dashed lines denote the Apollonius domains with colors matching their respective pursuers.
    }
    \label{fig:2v1-example.pdf}
\end{figure}

In any case, the agents' equilibrium strategy is to take a geodesic path to the equilibrium intercept point $I^*$ as descibed in \eqref{eq:two-v-one-solution}.
This is due to the fact that $I^* \in \partial \mathcal{A}_1 \cup \partial \mathcal{A}_2$ and all of the points on the boundary of an Apollonius domain are, by definition, points of simultaneous arrival via geodesic paths (i.e., see~\cite{Kovshov2000}).
Note that the state space for Proposition~\ref{prop:2v1} is restricted such that the geodesic distances from the pursuers to the evader satisfy~\eqref{eq:condThm3}.
This is because for $α > α_C$ the Apollonius domain $\mathcal{A}$ does not contain the equilibrium intercept point for the $\min$-$\max$ capture time game, and thus $\mathcal{E}$ may not necessarily be guaranteed to contain the equilibrium intercept point for the two-pursuer version of the game.
A rigorous proof for Proposition~\ref{prop:2v1}, as well as the solution for more general two-pursuer initial conditions is left for future work.

\subsection{Target Guarding}
\label{sec:Target_Guarding}

Another primary application for the Apollonius domain is to address target guarding games, in which an evader seeks to reach a target region while avoiding capture by one or more pursuers.
The evader wins if it can safely reach a point in the target region while the pursuer(s) wins if capture occurs before the evader can reach the target.
Let the target region be a subset of the surface of the sphere denoted as $\mathcal{T}$.

\begin{lemma}
    In the spherical target guarding game with one evader and one or more faster pursuers, the evader wins if $\mathcal{A} \cap \mathcal{T} \neq \varnothing$.
\end{lemma}

\begin{lemma}
    In the spherical target guarding game with one evader and one faster pursuer, the pursuer wins if
    \begin{equation}
        \label{eq:Kovshov-eq30} 
        \mathcal{A} \cap \mathcal{T} = \varnothing,
        \quad
        α \leq \tfrac{1 - μ}{1 + μ}π,
    \end{equation}
    by implementing the \textit{geodesic parallel strategy} (as defined in~\cite{Kovshov2000}), wherein the pursuer heads to the point on the Apollonius domain associated with the evader's current heading.
\end{lemma}
\begin{proof}
    From~\cite[Proposition 4.6]{Kovshov2000}, under the condition~\eqref{eq:Kovshov-eq30} if the pursuer uses the \textit{geodesic parallel strategy}, the Apollonius domain at the initial time is guaranteed to contain the Apollonius domain at all future times until capture for any evader's motion consisting of a polygonal line of geodesic arcs.
    Thus, if $\mathcal{A} \cap \mathcal{T} = \varnothing$ in the initial configuration, the condition is guaranteed to hold over the playout of the game.
    Moreover, from~\cite[Corollary to Proposition 4.7]{Kovshov2000} capture is guaranteed to occur within a time of $\frac{α}{1 - μ}$.
    Therefore, the evader's safe reachable set will never intersect the target region, capture is guaranteed in a finite time, and thus the evader will not be able to win.
\end{proof}

\section{Conclusion}
\label{sec:Conclusion}
This paper contains the solution to the pursuit-evasion differential game of $\min$-$\max$ capture time between one pursuer and one slower evader, all moving with simple motion on a sphere. In the special case where the two agents are on direct opposite ends of the sphere, the agents' equilibrium strategies are not well defined according to the conventional analysis; this configuration lies on a Dispersal Surface. A loss-rate analysis was carried out in this paper to determine a unique equilibrium control input for this configuration which involves the evader standing still for an infinitesimally small amount of time while the pursuer freely chooses any heading. The solution of the game provides a spherical analogue to the Apollonius circle, i.e., the Apollonius domain. It was found that the equilibrium intercept point lies outside the Apollonius domain in a portion of the state space (unlike in the planar case). Consequently, when the initial configuration lies within a specified parameter regime, the Apollonius domain exhibits similar properties to the Apollonius circle, i.e., the sphere can be considered flat. This is due to the geometry of a sphere, which is a closed domain, rather than an open domain of a Cartesian plane. Therefore, some caution must be taken when applying the concept of the Apollonius domain to solve spherical pursuit-evasion problems. Some example extensions and applications were discussed, including two cooperative pursuers and target guarding. Future work includes rigorously proving the two-pursuer solution, generalizing it, and obtaining evader reachability sets when the pursuer does not know the evader's heading. Establishing dominance regions for multi-pursuer-evader scenarios on spherical geometry remains to be a rich area of research.

\bibliographystyle{IEEEtran}
\bibliography{refs.bib}

\appendix[Relative $P$-$E$ Kinematics]

Here we derive the rate of change of angular distance $\alpha$, 
$\alpha \in (0,\pi)$ from \eqref{eq:nGCalpha} for $P$ and $E$ 
that obey \eqref{eq:evader}-\eqref{eq:pursuer}.

The time derivative of \eqref{eq:nGCalpha} yields
\begin{equation}
\dot{\vec{n}}_{GC} \sin \alpha + \vec{n}_{GC} \dot{\alpha} \cos \alpha = \tfrac{1}{R^2}\left(\dot{\vec{P}} \times  \vec{E} + 
\vec{P} \times  \dot{\vec{E}}\right). \label{eq:PdotEdot}
\end{equation}
Note that the magnitude $|\vec{n}_{GC}|=
\sqrt{\vec{n}_{GC}\cdot \vec{n}_{GC}}=1$, where ``$\cdot$" denotes the scalar product 
of vectors, therefore,
\begin{equation}
\frac{d|\vec{n}_{GC}|}{dt}=2 \dot{\vec{n}}_{GC}\vec{n}_{GC}=0
\Rightarrow \dot{\vec{n}}_{GC}\vec{n}_{GC}=0.
\end{equation}
Because of that, we apply the scalar 
product $\vec{n}_{GC}$ to both  sides of \eqref{eq:PdotEdot} and obtain  
\begin{equation}
\dot{\alpha} \cos \alpha = \tfrac{1}{R^2}(\dot{\vec{P}} \times  \vec{E}) \cdot \vec{n}_{GC} + 
\tfrac{1}{R^2} (\vec{P} \times  \dot{\vec{E}}) \cdot \vec{n}_{GC}.  \label{eq:dotalphacos}
\end{equation}

To find $\dot{\vec{P}}$ and $\dot{\vec{E}}$, let us introduce a reference frame composed of 
orthogonal unit vectors $\tfrac{1}{R}\vec{P}$, $\vec{n}_{GC}$ and $\vec{t}_P=\tfrac{1}{R}\vec{P} \times \vec{n}_{GC}$. 
In that reference frame, the $P$ and $E$ position vectors are
\begin{align}
 \vec{P} &= R \left(\tfrac{1}{R}\vec{P}\right) \\
 \vec{E} &=  R\cos \alpha\left(\tfrac{1}{R}\vec{P}\right)-
   R\sin \alpha \left(\tfrac{1}{R}\vec{P} \times \vec{n}_{GC}\right),
\end{align}
and the velocity vectors $\dot{\vec{P}}=v_P \vec{p}$ and $\dot{\vec{E}}=v_E \vec{e}$ are
\begin{align}
  \dot{\vec{P}} &= - v_P \cos u_P \left(\tfrac{1}{R}\vec{P} \times \vec{n}_{GC}\right) - v_P \sin u_P (\vec{n}_{GC}), \\
  \dot{\vec{E}} &= 
  \begin{multlined}[t]
      - v_E \sin \alpha \cos u_E \left(\tfrac{1}{R}\vec{P}\right)
      -v_E \sin u_E (\vec{n}_{GC}) \\
    - v_E \cos \alpha \cos u_E \left(\tfrac{1}{R}\vec{P} \times \vec{n}_{GC}\right),
  \end{multlined}
\end{align}
which also account for the relation between $\vec{t}_E$ and $\vec{t}_P$ 
(see~Fig.~\ref{Fig2}). 

Using vector identities, the two terms on the right side of 
\eqref{eq:dotalphacos} are
\begin{align}
(\dot{\vec{P}} \times  \vec{E}) \cdot \vec{n}_{GC} &= 
\begin{vmatrix}
0 & 1 &0 \\
0 & - v_P \sin u_P & - v_P \cos u_P \\
R\cos \alpha & 0 & -R\sin \alpha 
\end{vmatrix} \nonumber  \\
&=  - R v_P \cos u_P \cos \alpha 
\end{align}
\begin{align}
    (\vec{P} \times  \dot{\vec{E}}) \cdot \vec{n}_{GC}  &= 
    \begin{multlined}[t]
        \\ 
        \hspace{-7em} \begin{vmatrix}
            0 & 1 &0 \\
            R & 0  & 0 \\
            - v_E \sin \alpha \cos u_E & -v_E \sin u_E & - v_E \cos \alpha \cos u_E 
        \end{vmatrix} 
    \end{multlined}
\nonumber  \\
 &=  R v_E \cos u_E \cos \alpha 
\end{align}
which yields 
\begin{equation}
\dot{\alpha} \cos \alpha = \tfrac{1}{R} (v_E \cos u_E - v_P \cos u_P) \cos \alpha,
\end{equation}
and by matching the expressions on both sides of the equation, we can 
finally conclude that 
\begin{equation}
\dot{\alpha} = \tfrac{v_E}{R} \cos u_E - \tfrac{v_P}{R} \cos u_P.
\end{equation}
\end{document}